\documentclass[EJP]{ejpecp} 

\usepackage[vcentermath,enableskew]{youngtab}
\usepackage[procnames]{listings}
\usepackage{color}

\definecolor{keywords}{RGB}{255,0,90}
\definecolor{comments}{RGB}{0,0,113}
\definecolor{red}{RGB}{160,0,0}
\definecolor{green}{RGB}{0,150,0}

\SHORTTITLE{The Pillowcase Distribution and Near-Involutions} 

\TITLE{The Pillowcase Distribution and Near-Involutions\thanks{Supported by CONACYT-Mexico, the Mexican Secretariat of Public Education, and Princeton University.}} 

\AUTHORS{%
  Rodolfo R\'ios-Zertuche\footnote{Max-Planck-Institut f\"ur Mathematik, Bonn.
    \EMAIL{rrios@mpim-bonn.mpg.de}}}


\AMSSUBJ{20C32} 
\AMSSUBJSECONDARY{32G15; 62E20} 

\SUBMITTED{June 26, 2014} 
\ACCEPTED{December 14, 2014} 

\ARXIVID{1211.0735v2} 
\ARXIVPASSWORD{4haqj} 

\VOLUME{0}
\YEAR{2012}
\PAPERNUM{0}
\DOI{vVOL-PID}

\ABSTRACT{In the context of the Eskin-Okounkov approach to the calculation of the volumes of the different strata of the moduli space of quadratic differentials, the important ingredients are the pillowcase weight probability distribution on the space of Young diagrams, and the asymptotic study of characters of permutations that near-involutions. In this paper we present various new results for these objects. Our results give light to unforeseen difficulties in the general solution to the problem, and they simplify some of the previous proofs.}

\DeclareMathOperator{\trace}{tr}

\DeclareMathOperator{\var}{Var}
\DeclareMathOperator{\vol}{vol}
\DeclareMathOperator{\Aut}{Aut}

\newcommand{\C}{\mathbb{C}}
\newcommand{\R}{\mathbb{R}}
\newcommand{\Z}{\mathbb{Z}}

\newcommand{\w}{\mathsf{w}}
\newcommand{\p}{\mathbf{p}}
\newcommand{\s}{\mathbf{s}}
\newcommand{\g}{\mathbf{g}}
\newcommand{\f}{\mathbf{f}}

\newcommand{\vari}{\texttt}

\begin{document}



\abstract{In the context of the Eskin-Okounkov approach to the calculation of the volumes of the different strata of the moduli space of quadratic differentials, the important ingredients are the pillowcase weight probability distribution on the space of Young diagrams, and the asymptotic study of characters of permutations that near-involutions. In this paper we present various new results for these objects. Our results give light to unforeseen difficulties in the general solution to the problem, and they simplify some of the previous proofs.}

\section{Introduction}
\label{sec:intro}

A \emph{quadratic differential} on a Riemann surface $S$ is a function defined on the tangent bundle $\omega:TS\to\C$ that locally looks like $\omega(z)=f(z)(dz)^2$, for some meromorphic function $f$ with at most simple poles. There is an interesting correspondence between these objects and flat surfaces that makes them useful in the study of many dynamical systems, notably billiards; see \cite{zorich} for a good introduction to the subject.

We will say that two quadratic differentials $\omega$ and $\eta$ defined on surfaces $S$ and $T$ are equivalent if there is a holomorphic diffeomorphism $\varphi:S\to T$ such that $\varphi^*\eta=\omega$. The set of equivalence classes is known as the \emph{moduli space}. Its \emph{strata} $\mathcal M_\nu$ are the subsets corresponding to quadratic differentials with fixed numbers of zeros and poles. We encode this information in a partition $\nu=(\nu_1\geq\nu_2\geq\cdots\geq\nu_k\geq1)$, in which a simple pole will be represented by a part $\nu_i$ equal to 1, a marked point will be a 2, and a zero of degree $n$ will correspond to a part $\nu_i=n-2$. By Riemann-Roch, the genus $g$ of the underlying Riemann surface is determined by
\[|\nu|-2\ell(\nu)=\textrm{zeros}-\textrm{poles}=4g-4,\]
where $|\nu|=\sum_i\nu_i$ and $\ell(\nu)$ is the number $k$ of parts in $\nu$.

The strata $\mathcal M_\nu$ of the moduli space are known to be complex manifolds of dimension $2g-1+\ell(\nu)$ or $2g-2+\ell(\nu)$, depending on whether a global square-root of the quadratic differential exists or not. Since a scalar multiple of a quadratic differential is another quadratic differential with the same number of zeros and poles, the strata are cones.

Local coordinates in $\mathcal M_\nu$ were constructed in \cite{kontsevichlyapunovexponents}. The trick is to pass to the 2-fold cover $\tilde S$ where an Abelian differential $\eta$ (locally $\eta=f(z)\,dz$) exists that is a global square-root of the quadratic differential in question, $\omega=\eta^2$, and to take as coordinates the periods of $\eta$. The \emph{periods} are defined to be the values of integrals of $\eta$ along curves representing the generators of the homology group $H_1(\tilde S,P;\Z)$ of the surface $\tilde S$ relative to the set $P$ of zeros, poles, and marked points. Using these coordinates, one can pull Lebesgue measure from $\R^{\dim\mathcal M_\nu}$ to $\mathcal M_\nu$ to define a volume element.

Using those same periods and the Riemann bilinear relations, a quad\-rat\-ic differential is seen to induce the area of the surface. One can thus take the set $\mathcal M_\nu^1$ of quadratic differentials in $\mathcal M_\nu$ that induce area $\leq 1$. This subset was shown to have finite volume \cite{veech,masur}.

The problem of the computation of the volume of different strata of the moduli space of quadratic differentials, that is, of the volume of $\mathcal M^1_\nu$, gained relevance as connections with Siegel-Veech constants and Lyapunov exponents of the geodesic flow in Teichm\"uller space were discovered; see for example \cite{zorich,eskinmasurzorich}. 

The method developed by Eskin and Okounkov \cite{pillow} gives a way to compute the volumes. Their method leveraged the correspondence of quadratic differentials with flat surfaces to construct a lattice $F$ inside $\mathcal M_\nu$ consisting of tiled surfaces. These tiled surfaces are precisely the connected branched coverings of the pillowcase orbifold $\mathfrak P=(\R^2/\Z^2)/\{\pm1\}$ with branching data $\nu\cup(2,2,\dots)$ on one of the conic points of the pillowcase and $(2,2,\dots)$ on the three others. (A quadratic differential on these coverings is obtained by pulling back the differential $dz^2$ defined on the pillowcase.) They proved that
\[\vol \mathcal M^1_\nu=\lim_{d\to\infty}\frac{\textrm{$\#$ of connected degree $d$ coverings of $\mathfrak P$ in $\mathcal M_\nu$}} {d^{\dim \mathcal M_\nu-1}}.\]

They formed the generating function $Z_\nu$ of possibly disconnected coverings of the pillowcase in $\mathcal M_\nu$ graded by degree and weighted by \[\frac{1}{\Aut(S)|},\] where $\Aut(S)$ is the automorphism group of the covering $S$. (The weighting becomes asymptotically negligible as the degree increases, since most coverings have no nontrivial automorphisms; see \cite{branchedcoverings}.) This generating function is related to the generating function of connected coverings essentially by the exponential function; see \cite[Chapter 3]{wilf}.
They thus showed that in order to find the volume, one needs to understand the leading terms of the asymptotics as $q\to 1$, $|q|<1$, of the series
\begin{align}
\label{eq:coveringcount}
 Z_\nu(q)&=\sum_{d\geq 1}q^d\sum_{S\xrightarrow{d} \mathfrak P}\frac{1}{|\Aut(S)|}\\
\label{eq:partitionfunction}
&=\sum_\lambda q^{|\lambda|/2}\left(\frac{\dim\lambda}{|\lambda|!}\right)^2
\mathbf f_{(\nu,2,2,\dots)}(\lambda)\,\mathbf f_{(2,2,\dots)}(\lambda)^3
\end{align}
and that this series is a quasimodular form. The equivalence of \eqref{eq:coveringcount} and \eqref{eq:partitionfunction} is classical and is explained for example in \cite{gromovwittenhurwitz},\cite[Section 1.3.7]{mythesis}. Here,
\begin{itemize}
\item$\nu$ is a partition that indexes the stratum of the moduli space: the orders of the zeros of the quadratic differentials in the stratum are encoded as $\nu_i-2$ ($\nu_i=2$ correspond to marked points, and we do allow simple poles corresponding to the parts $\nu_i=1$),
\item the sum is taken over all partitions $\lambda=(\lambda_1,\lambda_2,\dots)$ of \emph{even} integers    $|\lambda|=\sum_i\lambda_i$,
\item $\dim \lambda$ denotes the dimension of the representation of the symmetric group corresponding to the partition $\lambda$,
\item $(2,2,\dots)$ means the partition with $|\lambda|/2$ parts equal to 2, and similarly 
\[(\nu,2,2,\dots)\] 
means the partition that contains all the parts of $\nu$ and is completed with twos to be of the same size as $\lambda$,
\item the numbers \[\f_{\eta}(\lambda)=|C_\eta|\frac{\chi^\lambda(\eta)}{\dim\lambda}\]
    are the central characters of the symmetric group algebra,
\item $|C_\eta|$ is the size of the conjugacy class of the symmetric group corresponding to permutations whose cycle type is given by the partition $\eta$, and
\item $\chi^\lambda(\eta)$ is the character of one such permutation in the irreducible representation of the symmetric group corresponding to the partition $\lambda$.
\end{itemize}

In practical terms, this means that in order to determine the volume of one of these strata one needs to compute the first few coefficients of the series \eqref{eq:partitionfunction}, from there deduce what polynomial in the Eisenstein series it corresponds to, and make use of the quasimodularity to determine the leading term in the asymptotics.

This is involved  because of the combinatoric calculation required in the first step, and it makes it impossible to conclude anything theoretically; it would be better to have a formula for the volumes. Our results go in this direction.

The Okounkov-Eskin method \cite{branchedcoverings,pillow} was further refined by Okounkov, Eskin, and Pandharipande \cite{eskinokounkovpandharipande}, but this refinement will not concern us here.

The alternative approach of Kontsevich \cite{kontsevich,kontsevichlyapunovexponents} and Kontsevich and Zorich \cite{kontsevichzorich} has been expanded and fruitfully exploited by Athreya, Eskin, and Zorich \cite{rightangledbilliards,jenkinsstrebeldifferentials}, who have found exact formulas for the volumes of several families of strata.

Let us state our results. The form of the series \eqref{eq:partitionfunction} inspires the following definition \cite{pillow}:

\begin{definition}[Pillowcase weight]\label{def:pillowcaseweight}
For a partition $\lambda$, its \emph{pillowcase weight} is defined by
\begin{equation*}
\w(\lambda)=\left(\frac{\dim \lambda}{|\lambda|!}\right)^2\f_{(2,2,\dots,2)}(\lambda)^4.
\end{equation*}
\end{definition}

The pillowcase weights induce a probability distribution in the space of Young diagrams: we introduce a parameter $q\in\C$, $|q|<1$, and we let the weight of the partition $\lambda$ be
\begin{equation}\label{eq:probabilitydistrib}
\frac{q^{|\lambda|/2}\w(\lambda)}{Z},
\end{equation}
where $Z=\sum_\mu q^{|\mu|/2}\w(\mu)$ is a normalization constant (independent of $\nu$).
Denoting the corresponding expectation by $\langle\cdot\rangle_{\w,q}$,
we can rewrite the series \eqref{eq:partitionfunction} as
\begin{equation}\label{eq:expectationofg}
Z_\nu(q)=Z\cdot \langle\g_\nu\rangle_{\w,q},
\end{equation}
where
\[\g_\nu(\lambda)=\frac{\f_{(\nu,2,2,\dots)}(\lambda)}{\f_{(2,2,\dots)}(\lambda)}.\]
The goal of this paper is to analyze the probability distribution \eqref{eq:probabilitydistrib} and to give some results about the computation of the expectations of $\g_\nu$.

Recall some definitions. The \emph{2-core} of a partition $\lambda$ is what remains of its Young diagram when one has removed as many 2-dominoes \[\yng(2)\qquad \yng(1,1)\]
as possible. 2-cores always have the shape $(k,k-1,\dots,1)$ of a staircase, for some $k$. We also associate to $\lambda$ a pair of partitions $(\alpha,\beta)$ that consitute its \emph{2-quotient}. They are defined as follows. Associate with $\lambda$ the sequence $\{z_i\}_{i\in\Z}$ of 0's and 1's such that $z_i$ is 1 if, and only if, there is $j>0$ such that $\lambda_j-j=i$. Then there is a unique $N\in\Z$ such that the sequence $\{z_{2(i+N)}\}_{i\in\Z}$ corresponds to a partition, and this is $\alpha$.  $\beta$ is the partition associated with the complementary sequence $\{z_{2(i-N)+1}\}_{i\in\Z}$.

For example, from the partition $(3,2,2)$, we can remove 2-dominoes as follows:
\[\young(\hfil\hfil\hfil,\hfil\times,\hfil\times)\to\young(\hfil\times\times,\hfil,\hfil)\to\young(\hfil,\times,\times)\to\young(\hfil)\]
We conclude that its 2-core is the partition $(1)$.
With the partition $(3,2,2)$, we associate the sequence of 0's and 1's given by
\[\dots,1,1,1,1,0,0,1,1,0,1,0,0,0,0,\dots\]
Separating the entries in even and odd positions, we obtain two sequences that correspond to the partitions $(1,1)$ and $(1)$, which together form the 2-quotient of $(3,2,2)$.

With the following definition, we characterize the support of $\w$.

\begin{definition}[Balanced partition]\label{def:balanced}
A partition is \emph{balanced} if either of the following equivalent statements is true:
\begin{itemize}
 \item $\w(\lambda)\neq 0$,
 \item the Young diagram of the partition $\lambda$ can be constructed by adjoining 2-do\-mi\-noes,
 \item the 2-core of $\lambda$ is empty,
 \item $\sum_i[(-1)^{\lambda_i-i+1}-(-1)^{-i+1}]=0$.
\end{itemize}
\end{definition}

For example, the partition $(4,2)$ is balanced, while $(3,2,2)$ is not --- its core is the partition (1) as shown above. The former can be constructed using 2-dominoes for example as follows:
\[\yng(1,1)\to\yng(2,2)\to\yng(4,2)\quad\textrm{or}\quad \yng(2)\to\yng(4)\to\yng(4,2)\]

Our first result is an interesting formula, proved in Section \ref{sec:hooksformula}:

\begin{proposition}\label{prop:formulawhooks}
For $\lambda$ balanced, we have
\[\w(\lambda)=\left(\frac{\prod\textrm{odd hook lengths of $\lambda$}}{\prod\textrm{even hook lengths of $\lambda$}}\right)^2.\]
\end{proposition}

For example, if $\lambda=(4,2)$, the hooklengths for each cell are:
\[\young(5421,21)\]
We obtain $\w(\lambda)=(5\cdot 1\cdot 1)^2/(4\cdot 2\cdot 2)^2=25/256$.

Using the formula in Proposition \ref{prop:formulawhooks}, we obtain the following estimate through a variational argument; see Section \ref{sec:variational} for the proof. 
\begin{proposition}\label{prop:closetodiagonal}
Let $\varepsilon>0$, and let $S_{\varepsilon,n}$ be the set of balanced partitions $\lambda$ of $n$ whose 2-quotients $(\alpha,\beta)$ satisfy $\|L_\alpha-L_\beta\|>\varepsilon$, where $\|\cdot\|$ denotes a Sobolev norm and
$L_\alpha$ and $L_\beta$ are the rescaled contours of the partitions $\alpha$ and $\beta$ (see Definitions \ref{def:contour} and \ref{def:sobolev}). Then the $\w$-probability of $S_{\varepsilon,n}$ is asymptotically of order $O\left(e^{-K\sqrt n}\right)$ as $n\to\infty$, for some $K>0$ that depends on $\varepsilon$. In fact $K\to0$ as $\varepsilon\searrow0$.
\end{proposition}

Whence we conclude that the measure is concentrated on the set of partitions that are close to having identical 2-quotient components $\alpha$ and $\beta$, that is, close to being constructed solely by $2\times2$ blocks:
\[\yng(2,2),\]
and that it must be nearly uniform there. We also expect this to be true since, for large partitions, the odd hook lengths should almost cancel out the even hook lengths in the formula of Proposition \ref{prop:formulawhooks}.

Although we are unable to extract the uniformity from the result in Proposition \ref{prop:closetodiagonal}, we can still prove it through the analysis of the expression for the $n$-point function obtained in \cite{pillow}. Hence we have the following result, proved in Section \ref{sec:limitshape}.

\begin{proposition}\label{prop:limitshape}
The probability distribution induced by the weights $\w(\lambda)$ induces a limit shape that coincides with the one for the uniform probability distribution, namely, it is the curve
\[e^{-\pi x/\sqrt 6}+e^{-\pi y/\sqrt6}=1.\]
\end{proposition}

This means that most of the time, partitions $\lambda$ sampled with that probability distribution will be very close to that curve. Here is a simulation:
\begin{figure}[htbp]
\begin{center}
\includegraphics[scale=1.0]{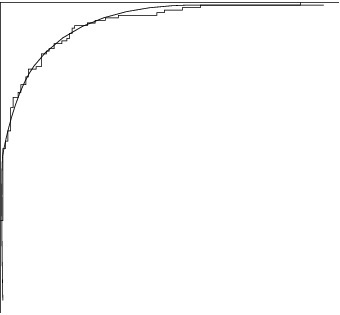}
\end{center}
\end{figure}

We were eager to find this result, because it determines immediately the asymptotic value of the expectations of the shifted power functions, and the function $\g_\nu$ turns out to be a polynomial in the shifted power functions \cite{pillow}. From there, we expected to find a Wick-type theorem that would further simplify our computation of the volumes. Instead we did some computations that strongly suggest that \emph{for the distribution given by \eqref{eq:probabilitydistrib}, the convergence to the limit shape is not normal. In other words, that there is no central limit theorem for the coupled pillowcase weights \eqref{eq:probabilitydistrib} as $q\to 1$, and that there is no Wick-type theorem associated to the distribution they induce.
} We discuss this in Section \ref{sec:noclt}.
\begin{remark}
It may be that the distribution induced by the ``uncoupled'' weights $\w(\lambda)$ (instead of the weights given in \eqref{eq:probabilitydistrib}) has normal convergence to the limit shape, and there may be a way to exploit this in order to compute the volumes. The strategy that our computations would rule out is the use of Wick's theorem together with an application of the result of Proposition \ref{prop:formulaforg} below to produce the volumes, if we knew only the asymptotic values of the expectations $\langle\p_i\p_j\rangle_{\w,q}$ of products of two shifted-power functions. A similar strategy does work in the case of the volumes of strata of moduli spaces of Abelian differentials \cite{branchedcoverings}.
\end{remark}

We also find the following structural formula, which will be proved in Section \ref{sec:nearinvolutions} and commented in Remarks \ref{rmk:quasimodularity} and \ref{rmk:matrix}.

\begin{proposition}\label{prop:formulaforg}
\begin{equation}
\g_\nu(\lambda)=\frac{2^{|\nu|/2}}{\mathfrak z(\nu)}\sum_\mu (-1)^{o_\mu/2}\chi^\mu(\nu)\s_a(\alpha)\s_b(\beta),
\end{equation}
where
\begin{itemize}
\item the sum is taken over all balanced partitions $\mu$ (see Definition \ref{def:balanced}) of size $|\mu|=|\nu|$ whose Young diagram is completely contained\footnote{Strictly speaking, the requirement that the diagram of $\mu$ be contained in the diagram of $\lambda$ is superfluous because at least one of the shifted Schur functions will vanish in case it is not.} inside the Young diagram of $\lambda$,
\item $o_\mu$ is the number of odd parts $\mu_i$ of $\mu$,
\item $\s_\mu$
are the shifted Schur functions introduced in Definition \ref{def:shiftedschur} below,
\item $\mathfrak z(\nu)=\prod_{n=1}^\infty n^{m_n}m_n!$
(where $m_n$ is the number of parts $\nu_i$ equal to $n$)
is the cardinality of the centralizer of an element of cycle type $\nu$ in the symmetric group, and
\item the pairs of partitions $(\alpha,\beta)$ and $(a,b)$ are the 2-quotients of $\lambda$ and $\mu$, respectively.
\end{itemize}
\end{proposition}

In the final Section \ref{sec:application}, we discuss the initial applications of these results to the calculation the expectations \eqref{eq:expectationofg}.

This paper elaborates on some of the results of the author's PhD thesis \cite{mythesis}, where a more detailed account of most of the proofs can be found.

\section{Formula in terms of hooks}
\label{sec:hooksformula}

\begin{proof}[Proof of Proposition \ref{prop:formulawhooks}]
Let $\lambda$ be a balanced partition. Note first that
\begin{equation}\label{eq:sizeofconjugacyclass}
|C_{(2^{|\lambda|/2})}|=\frac{|\lambda|!}{2^{|\lambda|/2}(|\lambda|/2)!},
\end{equation}
and recall the classical result that \cite{pillow},\cite[Section 2.2.1]{mythesis}
\[|\chi^\lambda(2,2,\dots,2)|=\binom{|\lambda|/2}{|\alpha|}\dim\alpha\dim\beta,\]
where the pair of partitions $(\alpha,\beta)$ is the 2-quotient of $\lambda$. Since $\lambda$ is balanced, exactly half of its hook lengths are even, and the hook lengths of $\alpha$ and $\beta$ are in one-to-one correspondence with the even hook lengths of $\lambda$ divided by 2.
(See \cite[Section 2.1]{mythesis}.)
Using this equality and the hook formula of Frame-Robinson-Thrall, namely,
\begin{equation}\label{eq:hookformula}
\dim\lambda=\frac{|\lambda|!}{\prod_{\square\in\lambda}h_\square},
\end{equation}
where $h_\square$ denotes the length of the hook of the cell $\square\in\lambda$,
 we get
\[\dim\alpha\dim\beta=\frac{|\alpha|!|\beta|!}{\prod_{\{\square\in\lambda: 2|h_\square\}} (h_\square/2)}\]
where the product is taken over those cells whose hook lengths are even. In the denominator we have the product of the halves of the even hook lengths of $\lambda$. In other words,
\begin{equation}\label{eq:involutionformula}
|\chi^\lambda(2,2,\dots,2)|=\frac{2^{|\lambda|/2}(|\lambda|/2)!}{\prod\textrm{even hook lengths of $\lambda$}}.
\end{equation}
Since
\begin{equation*}
\frac{\prod\textrm{hook lengths of $\lambda$}}{\left(\prod\textrm{even hook lengths of $\lambda$}\right)^2} =\frac{\prod\textrm{odd hook lengths of $\lambda$}}{\prod\textrm{even hook lengths of $\lambda$}},
\end{equation*}
plugging \eqref{eq:sizeofconjugacyclass}, \eqref{eq:hookformula} and \eqref{eq:involutionformula} into the expression given in Definition \ref{def:pillowcaseweight} for $\w(\lambda)$, we get the desired formula.
\end{proof}

\section{Concentration near the diagonal}
\label{sec:variational}
The goal of this section is to prove Proposition \ref{prop:closetodiagonal}. We first collect some definitions and basic results.

\begin{definition}[Rescaled rim function]
Let $\lambda$ be a partition of $n=|\lambda|$. We contract its diagram until it has area 1, rescaling by $1/\sqrt n$, and we associate to it the \emph{rescaled rim function} of $\lambda$, namely, the non-increasing function $F:\R_+\to\R_+$ given by
\[F(x)=\frac{1}{\sqrt n}\#\{\textrm{parts of $\lambda$ of size $\leq \lceil\sqrt n\,x\rceil$\}}.\]
\end{definition}

\begin{definition}[Approximate hook]
For an increasing function $F:\R_+\to\R_+$, define
\[F^{-1}(y)=\inf\{x:F(x)\leq y\},\]
and
\[h_F(x,y)=F(x)+F^{-1}(y)-x-y,\]
its \emph{approximate hook} at $(x,y)$.
\end{definition}

\begin{lemma}[Kerov-Vershik \cite{vershikkerov1985}]\label{lem:integralapproximation}
Let $\lambda$ be a partition and let $\square\in \lambda$ be a cell in its Young diagram, and denote by $F$ its rescaled rim function. Let $R_\square\subseteq\{(x,y)\in \R^2:0\leq y\leq F(x)\}$ be the rectangular domain (of area $1/n$) corresponding to $\square$. Let $h_\square$ denote the hook length of the cell $\square$. Then
\[\log h_\square=n \int_{R_\square}\log\left(\sqrt nh_F(x,y)\right)dx\,dy+c(h_\square),\]
where
\[c(x)=\frac12\sum_{k=1}^\infty\frac{1}{k(k+1)(2k+1)x^{2k}}.\]
\end{lemma}

\begin{lemma}\label{lem:orderofc}
\[\sum_{\square\in\lambda}(-1)^{h_\square}c(h_\square)=O\left(\sqrt{|\lambda|} \right) \quad as \quad|\lambda|\to\infty.\]
\end{lemma}
\begin{proof}[Proof of Lemma \ref{lem:orderofc}]
We have
\begin{equation}\label{eq:boundforc}
\left|\sum_{\square\in\lambda}(-1)^{h_\square}c(h_\square)\right|\leq
\sum_{\{\square\in\lambda:\textrm{$h_\square$ odd}\}}c(h_\square).
\end{equation}
Since $c$ is decreasing and $c(h)\to\infty$ as $h\to\infty$, the worst case scenario is the case in which we maximize the number of \emph{small} odd hook lengths $h_\square$. This happens in the case of the staircase partition $(\ell,\ell-1,
\dots,2,1)$, and in this case the right hand side of \eqref{eq:boundforc} is of order $O\left(\sqrt{|\lambda|}\right)$.
\end{proof}

\begin{definition}[Rescaled contour of a partition]\label{def:contour}
The \emph{rescaled contour} of the Young diagram of a partition $\lambda$ is the function $L_\lambda:\R\to\R_+$ whose graph is the union of the graphs of the rescaled rim function $F$ of $\lambda$ and of $F^{-1}$, after the change of variables $x'=(y-x)/\sqrt2$, $y'=(x+y)/\sqrt2$.

In other words, we rotate the rescaled diagram of $\lambda$ so that its sides become aligned with the graph of $x\mapsto |x|$; then $L_\lambda$ is the function that describes the rim of the rotated diagram, and it equals $|x|$ outside the diagram.
\end{definition}

\begin{definition}[Sobolev norm \cite{vershikkerov1985}]\label{def:sobolev}
Define by
\[\|f\|=\int_{\R^2}\left(\frac{f(s)-f(t)}{s-t}\right)^2ds\,dt.\]
the \emph{Sobolev norm} of $f$.
\end{definition}

\begin{proposition}\label{prop:firstestimate}
As $|\lambda|\to \infty$ (and restricting to $\lambda$ balanced),
\[\w(\lambda)=\exp\left(-\frac{|\lambda|}2\|\Delta\|^2+O\left(\sqrt{|\lambda|}\right) \right),\]
where $\Delta=L_\alpha-L_\beta$ is the difference of the rescaled contours of the components $\alpha$ and $\beta$ of the 2-quotient of $\lambda$.
\end{proposition}
\begin{proof}[Proof of Proposition \ref{prop:firstestimate}]
Consider the diagram of a balanced partition $\lambda$ rescaled by $1/\sqrt n$.
Denote by $O$ and $E$ the domains inside this diagram corresponding to the cells whose hooks are of odd and even length, respectively. It follows from Lemma \ref{lem:integralapproximation} that
\[\log\w(\lambda)=2n\left(\int_O-\int_E\right)\log\left(\sqrt n\,h_F(x,y)\right)dx\,dy+ \sum_{\square\in\lambda}(-1)^{h_\square}c(\lambda).\]
In a balanced partition, the area of $O$ is the same as the area of $E$, so the contributions of $\log\sqrt n$ cancel out. So $\sqrt n$ can be removed from the integrand, and from Lemma \ref{lem:orderofc} we have
\[\log\w(\lambda)=2n\left(\int_O-\int_E\right)\log\left(h_F(x,y)\right)dx\,dy+ O\left(\sqrt{|\lambda|}\right).\]
We change variables to $s$ and $t$ such that $x=\frac1{\sqrt 2}(L_\lambda(s)-s)$ and $y=\frac1{\sqrt 2}(L_\lambda(t)+t)$, and the first term above becomes
\[n\left(\int_{\tilde O}-\int_{\tilde E}\right)\left(\log\sqrt2(s-t)\right)(1-L_\lambda'(s))(1+L_\lambda'(t))\,ds\,dt,\]
where $\tilde O$ and $\tilde E$ are the images of $O$ and $E$ under the change of variables.

Let $(\alpha,\beta)$ be a pair of partitions giving the 2-quotient of $\lambda$. Since the Maya diagram of $\lambda$ is a sequence that can be obtained from the Maya diagrams of $\alpha$ and $\beta$ by placing their elements alternatingly, we have
\begin{equation}\label{eq:contourderivative}
L_\lambda'(s)=\left\{\begin{array}{ll}
L_\alpha'\left(s-\frac{k}{2\sqrt{2n}}\right),& \textrm{for $s\in\frac1{\sqrt{2n}}(k-1,k),k\in 2\Z$,}\\
L_\beta'\left(s-\frac{k}{2\sqrt{2n}}\right),& \textrm{for $s\in\frac1{\sqrt{2n}}(k,k+1),k\in 2\Z$.}
\end{array}\right.
\end{equation}

Again referring to the Maya diagram of $\lambda$, note that a hook of length $k$ corresponds to a pair of slots $k$ units apart such that the slot on the left has a pebble and the slot on the right does not \cite[Excercise 7.59]{stanley}. In particular, it follows that the integral over $\tilde E$ involves the interaction of each component of the 2 quotient with itself, so this term becomes
\begin{multline}\label{eq:integraloverE}
-n\int_{t<s}\left(\log\sqrt2(s-t)\right)\left(\left(1-L'_\alpha(s)\right) \left(1+L'_\alpha(t)\right)
\right.\\ \left.
+\left(1-L'_\beta(s)\right)\left(1+L'_\beta(t)\right)\right)ds\,dt.
\end{multline}
Similarly, the integral over $\tilde O$ involves interactions of the components of the 2-quo\-tient with each other, and becomes
\begin{multline}\label{eq:integraloverO}
n\int_{t<s}\left(\log\sqrt2(s-t)\right)\left(\left(1-L'_\alpha(s)\right) \left(1+L'_\beta(t)\right)
\right.\\ \left.
+\left(1-L'_\beta(s)\right)\left(1+L'_\alpha(t)\right)\right)ds\,dt.
\end{multline}
Note the implicit change of coordinates here, which would handle the translations by $k/\left(2\sqrt{2n}\right)$ in formula \eqref{eq:contourderivative}.

Adding \eqref{eq:integraloverE} and \eqref{eq:integraloverO}, and simplifying, we obtain
\[n\int_{\R^2}\left(\log\sqrt2|s-t|\right)\Delta'(s)\Delta'(t)\,ds\,dt,\]
where $\Delta(s)=L_\alpha(s)-L_\beta(s)$, from where the statement of the proposition follows after integrating by parts twice.
\end{proof}

\begin{lemma}\label{lem:meinardus}
As $n\to\infty$, we have the asymptotic behavior
\[Z_n=\sum_\lambda\w(\lambda)\sim\frac{1}{2^{1/8}3^{3/8}}\frac{e^{\pi\sqrt{n/6}}}
{n^{7/8}},\]
where the sum is taken over all $\lambda$ of size $|\lambda|=n$.
\end{lemma}
\begin{proof}[Proof of Lemma \ref{lem:meinardus}]
It was shown in Eskin-Okounkov \cite{pillow} that
\[\sum_\lambda\w(\lambda)q^{|\lambda|}=\prod_{i=1}^\infty\left(1-q^{2i}\right)^{-1/2}.\]
We have thus the following parameters for the Theorem of Meinardus, as presented by Andrews \cite[Chapter 6]{andrews}: $q$ in the book is $q^2$ here, $a_n=\frac12$, $D(s)=\zeta(s)/2$, $\alpha=1$, $A=\frac12$, $\kappa=-\frac78$, $D(0)=-\frac14$, $D'(0)=\frac14\log2\pi$. The statement of the lemma follows.
\end{proof}

\begin{proof}[Proof of Proposition \ref{prop:closetodiagonal}]
We can bound the cardinality of $S_{\varepsilon,n}$ by the total number $p(n)$ of partitions of $n$, for which the asymptotics
\begin{equation}\label{eq:asymptoticp}
p(n)\sim\frac1{4n\sqrt3}\exp\left(\pi\sqrt{\frac23n}\right)
\end{equation}
is well known \cite[Theorem 6.3]{andrews}.
The $\w$ probability of $\lambda$ is given by
\[\frac{\w(\lambda)}{\sum_{|\mu|=n}\w(\mu)},\]
 The statement of the proposition follows from Proposition \ref{prop:firstestimate}, Lemma \ref{lem:meinardus}, and the asymptotic relation \eqref{eq:asymptoticp}.
\end{proof}

\section{The limit shape}
\label{sec:limitshape}
Our goal in this section is to prove Proposition \ref{prop:limitshape}. We first need a few results and definitions.

Recall that the rescaled contour $L_\lambda$ of $\lambda$ was introduced in Definition \ref{def:contour}.

\begin{definition}[Limit shape]
Let $\nu$ be a measure on the set of Young diagrams, such that $\nu(\{\lambda:|\lambda|=n\})=1$ for all $n>0$.
A continuous function $\Omega:\R\to\R_+$ such that
\[\int_\R\left[\Omega(x)-|x|\right]\,dx=1\]
is the \emph{limit shape} induced by $\nu$ if there is some function $e:\R_+\to\R_+$ such that $e(x)\to0$ as $x\to0$ and for all $\varepsilon>0$ there is some $N\gg0$ such that, for all $n>N$,
\[\nu(\{\lambda:|\lambda|=n, \|L_\lambda-\Omega\|_\infty<e(\varepsilon)\})>1-\varepsilon.\]
\end{definition}

\begin{definition}[Moments of the limit shape]
For a limit shape $\Omega$, we define its \emph{moments}
\[\mu_k=\int_\R x^k(\Omega(x)-|x|)\,dx.\]
\end{definition}

\begin{definition}[Shifted power functions]
Let
\[\p_k(\lambda)=\left(1-\frac1{2^k}\right)\zeta(-k)+\sum_{i}\left[\left(\lambda_i-i+\tfrac12\right)^k- \left(-i+\tfrac12\right)^k\right].\]
The \emph{shifted power function} $\p_\mu$ indexed by the partition $\mu$ is defined as
\[\p_\mu(\lambda)=\prod_k \p_{\mu_k}(\lambda).\]
\end{definition}

\begin{remark}\label{rmk:moments}
When a distribution $\nu$ on the set of Young diagrams has a limit shape,
\[\lim_{n\to\infty}\langle\p_k\rangle_{\nu,n}=k\mu_{k-1}\]
where $\langle\cdot\rangle_{\nu,n}$ denotes the mean among all partitions of size $|\lambda|=n$. (See for example \cite{ivanovolshanski}, \cite[Section 3.1]{mythesis}.)
\end{remark}

Most of our derivation will come from Theorem \ref{thm:eskok} below. In order to state it, we need some more definitions.
\begin{definition}[Jacobi theta function]
The \emph{Jacobi theta function} is given by
\[\vartheta(x)=\vartheta(x,q)=\left(x^{1/2}-x^{-1/2}\right)\prod_{i=1}^\infty \frac{(1-q^ix)(1-q^i/x)}{(1-q^i)^2}.\]
\end{definition}
\begin{remark}
Recall the definition \cite{mumford} of the classical Jacobi theta functions $\theta_{jk}$, $j,k=0,1$:
\[\vartheta_{jk}(e^u,e^{-h})=\sum_{n\in\Z}\exp\left(-\left(n+\frac j2\right)^2\frac h2+ \left(n+\frac j2\right)\left(u+\pi i k\right)\right).\]
In terms of these and as a consequence of the Jacobi triple product formula,
\begin{equation}\label{eq:classictheta}
\vartheta(x,q)=-\eta^{-3}(q)\vartheta_{11}(x,q),
\end{equation}
where $x=e^u$, $q=e^{-h}$, and $\eta$ stands for the Dedekind eta function,
\[\eta(q)=q^{1/24}\prod_{n=1}^\infty(1-q^n).\]
\end{remark}

\begin{definition}
The \emph{$n$-point function} is the generating function
\[F(x_1,\dots,x_n)=\frac1Z\sum_{\lambda}q^{|\lambda|/2}\w(\lambda) \left(\sum_{i_1} x_1^{\lambda_{i_1}-i_1+\tfrac12}\right)\cdots\left(\sum_{i_n} x_n^{\lambda_{i_n}+i_n+\tfrac12}\right)\!.\]
\end{definition}

\begin{remark}
Since \cite{blochokounkov,pillow}
\[\sum_i e^{u(\lambda_i-i+\tfrac12)}=\frac1u+\sum_k\frac{u^k}{k!} \p_k(\lambda),\]
we have
\[F(e^{u_1},\dots,e^{u_n})=\sum_{i_1,\dots,i_n=1}^\infty
\frac{u_1^{i_1}\cdots u_n^{i_n}}{i_1!\cdots i_n!}
\langle \p_{i_1}\cdots\p_{i_n}\rangle_{\w,q}+\textrm{irrelevant terms}.
\]
(The terms missing in this formula are irrelevant in the sense that they do not encode information of interest to us. They correspond to monomials in the variables $u_i$ that are different from the ones in the sum.)
Hence, in view of Remark \ref{rmk:moments}, a statement about $F$ should be interpreted as a statement about the moments of the limit shape.
\end{remark}

\begin{theorem}[Eskin-Okounkov \protect{\cite[Theorem 5]{pillow}}]
\label{thm:eskok}
We have
\begin{equation*}
F(x_1,\dots,x_n)=\prod\frac{1}{\vartheta(x_i)}
[y^0_1\cdots y_n^0]\prod_{i<j}\frac{\vartheta(y_i/y_j)\vartheta(x_i y_i/x_jy_j)}{\vartheta(x_iy_i/y_j)\vartheta(y_i/x_jy_j)}
\prod_i\sqrt{
\frac{\vartheta(-y_i)\vartheta(x_iy_i)}{\vartheta(y_i)\vartheta(-x_iy_i)}},
\end{equation*}
where the brackets $[\cdot]$ indicate the operation of taking the coefficient of the indicated monomial, and the series expansion is performed in the domain
\[|y_n/q|>|x_1y_1|>|y_1|>\cdots>|x_ny_n|>|y_n|>1.\]
\end{theorem}

\begin{proposition}\label{prop:asymptoticsnpoint}
The highest degree term in the asymptotics of the $n$-point function is
\[F\left(e^{hu_1},e^{hu_2},\dots,e^{hu_n}\right)=\left(\prod_{i=1}^n \frac1{2h\sin\left(\pi u_i/2 \right)}\right)\left(1+O(h)\right)\]
as $h\to+0$.
\end{proposition}

For the proof of Proposition \ref{prop:asymptoticsnpoint} we need two lemmas.
We will use the notation $\approx$ to mean ``up to exponentially small terms.''

\begin{lemma}[\protect{\cite[Proposition 4.1]{branchedcoverings}}]\label{lem:theta1}
We have
\[\frac{\vartheta(e^{u},e^{-h})}{\vartheta'(0,e^{-h})}\approx
h\frac{\sin(\pi u/h)}{\pi}\exp\left(\frac{u^2}{2h}\right)\]
as $h\to+0$, uniformly in $u$. This asymptotic relation can be differentiated any number of times.
\end{lemma}

\begin{lemma}\label{lem:theta2}
We have the following approximations for $\vartheta$ at $e^u=1$ and $e^u=-1$, respectively:
\begin{align*}
\vartheta(e^u,e^{-h})&\approx i\eta^{-3}(q)\sqrt{\frac{2\pi}h}\exp\left(\pi i \left(\left[\frac{u}{2\pi i}-\tfrac12\right]+\tfrac12\right)-\frac{2\pi^2}h\left\{\frac{u}{2\pi i}-\tfrac12\right\}^2\right),\\
\vartheta(-e^u,e^{-h})&\approx i\eta^{-3}(q)\sqrt{\frac{2\pi}h}\exp\left(-\frac{2\pi^2}h\left\{\frac{u}{2\pi i}-\tfrac12\right\}^2\right),
\end{align*}
as $h\to0$, where $q=e^{-h}$, $[x]$ stands for the integer closest to $x$,  and $\{x\}=x-[x]$.
\end{lemma}
\begin{proof}[Proof of Lemma \ref{lem:theta2}]
This is a straightforward application of the modular transformation $h\mapsto -1/h$ in the expression \eqref{eq:classictheta}. We also use the identity \[\vartheta_{11}(-e^{u},q)=-\vartheta_{10}(e^u,q).\qedhere\]
\end{proof}

\begin{proof}[Proof of Proposition \ref{prop:asymptoticsnpoint}]
We start with the expression from Theorem \ref{thm:eskok}.
We first need to understand the asymptotic behavior of the coefficient of $y_1^0\cdots y_n^0$. We will approach this as an integral
\[\oint_{|y_1|=c}\frac{dy_1}{y_1}\cdots \oint_{|y_n|=c}\frac{dy_n}{y_n}\,Q,\]
where $Q$ is the quotient of theta functions to the right of $[y_1^0\cdots y_n^0]$ in Theorem \ref{thm:eskok}, including the part with the square root.
The approximation of Lemma \ref{lem:theta2} implies that if $x_i=e^{u_ih}$, $y_i=e^{v_i}$, and $q=e^{-h}$, then $Q$ is the exponential of a sum of many terms. Some of them vanish as $h\to0$; for the rest, we have to take l'H\^opital's rule. In the end, we see that $Q$ tends to
\[\exp\left(\pi i\sum_i u_i\left(\left\{ \frac{v_i}{2\pi i}-\tfrac12\right\}-\left\{ \frac{v_i}{2\pi i}\right\}\right)\right). \]
We take the integral for each $v_i$ on the segment  $[0,2\pi i]\subset \C$, and we get
\[\prod_i \pi\cos(\pi u_i).\]

Finally we apply the approximation of Lemma \ref{lem:theta1} to the product $\prod_i 1/\vartheta(x_i)$ in front of the expression of Theorem \ref{thm:eskok}.
\end{proof}

\begin{corollary}\label{coro:multiplicative}
The expectations of the shifted power functions are asymptotically multiplicative, namely
\[\langle\p_{i_1}\cdots\p_{i_k}\rangle_{\w,q}= \langle\p_{i_1}\rangle_{\w,q}\cdots\langle\p_{i_k}\rangle_{\w,q}
+o\left(h^{-\sum_k(i_k+1)}\right)\quad \textrm{as $h\to+0$}.\]
\end{corollary}
\begin{proof}[Proof of Corollary \ref{coro:multiplicative}]
This follows from Proposition \ref{prop:asymptoticsnpoint} and from the estimate \cite[Section 3.3.5]{pillow}
\[\langle\p_\rho\rangle_{\w,q}=O\left(h^{-|\rho|-\ell(\rho)}\right),\]
where $\ell(\rho)$ denotes the number of parts of the partition $\rho$.
\end{proof}

\begin{proof}[Proof of Proposition \ref{prop:limitshape}]
In order to determine that the limit shape should be equal to the one induced by the uniform distribution, we will show that, in the $q\to1$ limit, the asymptotic behavior of the expectations $\langle\p_k\rangle_{\w,q}$ coincides with $\mu_{k-1}^{\textrm{unif}}/k$, where $\mu_k^{\textrm{unif}}$  are the moments of that limit shape (compare with Remark \ref{rmk:moments}). Because these moments do not grow very fast, they determine the shape uniquely. A sufficient condition \cite[Chapter 30]{billingsley} for this to be the case is that
\[\limsup_{k\to\infty} \frac1{2k}(\mu_{2k})^{1/2k}<\infty,\]
which is true in this instance.

In the case of the $q$-coupled uniform distribution, in which the weight of the partition $\lambda$ is $q^{|\lambda|}$, the corresponding $n$-point function (which is also the exponential generating function of the numbers $\langle\p_k\rangle_{\textrm{unif},q}$) is
\[
F_\textrm{unif}(x)=\frac1u+\sum_{j\geq0}\frac{u^j}{j!}\left\langle\p_j \right\rangle_{\textrm{unif},q}=\frac1{\vartheta(x)},
\]
where $x=e^u$.The proof is similar to that of Theorem \ref{thm:eskok}, noting that (in the notations of \cite{pillow})
\[F_\textrm{unif}(x)=\frac{1}{\sum_\lambda q^{|\lambda|}}[y^0]\trace q^H\psi(xy)\psi^*(y).\]
Asymptotically, this behaves as
\[F_\textrm{unif}(e^{hu})= \frac1{\vartheta(e^{hu})} \approx\frac{\pi}{h\sin\pi u}=\frac1h\left(\frac1u+\sum_{j}\frac{u^j}{j!}\mu_j^{\textrm{unif}}\right),
\]
where $q=e^{-h}$, $h\to+0$, and $\approx$ means ``up to exponentially small terms.''
These asymptotics are obtained through an application of Lemma \ref{lem:theta1}.

From Proposition \ref{prop:asymptoticsnpoint} we see that the $1$-point function $F$ for the pillowcase distribution is (up to a rescaling) asymptotically equivalent to $F_\textrm{unif}$, whence the moments of the limit shape are also the same, $\mu_j=\mu_j^\textrm{unif}$.

On the other hand, from the multiplicativity property of Corollary \ref{coro:multiplicative} it follows that for any partition $\rho$, the variance
\[\var \p_\rho=\langle\p_\rho^2\rangle_{\w,q}-\langle\p_\rho\rangle_{\w,q}^2\to 0\]
very quickly as $q\to1$. This means that the probability must be concentrated at one point. In other words, our candidate is indeed the limit shape.
\end{proof}

\section{No central limit theorem}
\label{sec:noclt}
In this section, we want to discuss some examples we computed that suggest that for the distribution given by \eqref{eq:probabilitydistrib}, the convergence to the limit shape is not normal.

\begin{remark}
These computations rely heavily on a computer program (presented in Appendix \ref{apx:program}). While we have no reason to doubt its correctness, we also currently have no way to validate it beyond the highest-weight terms, which do coincide with the results of Propositions \ref{prop:closetodiagonal} and \ref{prop:limitshape}.
\end{remark}

If the probability measure were asymptotically gaussian, a Wick-type theorem would hold for the shifted power functions $\p_\mu$. Namely, we would have an identity of the type
\[\langle abcd\rangle_{\w,q}=\langle ab\rangle_{\w,q}\langle cd\rangle_{\w,q} +\langle ac\rangle_{\w,q}\langle bd\rangle_{\w,q} +\langle ad\rangle_{\w,q}\langle bc\rangle_{\w,q}\]
for all functions $a$, $b$, $c$, and $d$ with vanishing $(\w,q)$-mean and contained in the algebra generated by the functions $\p_k$. In particular, this means that we would need to have
\begin{equation}\label{eq:needforp1}
\left\langle \left(\p_1-\langle\p_1\rangle_{\w,q}\right)^4\right\rangle_{\w,q}
=3\left\langle \left(\p_1-\langle\p_1\rangle_{\w,q}\right)^2\right\rangle_{\w,q}^2
\end{equation}
and
\begin{multline}\label{eq:needforp1p3}
\left\langle \left(\p_1-\langle\p_1\rangle_{\w,q}\right)^3
\left(\p_3-\langle\p_3\rangle_{\w,q}\right)
\right\rangle_{\w,q}\\
=3\left\langle \left(\p_1-\langle\p_1\rangle_{\w,q}\right)^2\right\rangle_{\w,q}
\left\langle \left(\p_1-\langle\p_1\rangle_{\w,q}\right)\left(\p_3-\langle\p_3\rangle_{\w,q}\right)\right\rangle_{\w,q}.
\end{multline}
To check whether this was the case, we computed the first few terms of the following series, and then we used the quasimodularity property proved by Eskin-Okounkov \cite[Theorem 1]{pillow} with the methodology detailed in Appendix \ref{apx:program} to get:
\begin{align*}
\langle\p_1\rangle_{\w,q}&=\frac{\pi ^2}{24 h^2}+\frac{1}{4 h}+\mathrm{e.s.t.} \\
\langle\p_1^2\rangle_{\w,q}&
 =\frac{\pi ^4}{576 h^4}-\frac{\pi ^2}{16 h^3}-\frac{3}{16 h^2}+\mathrm{e.s.t.} \\
\langle\p_1^3\rangle_{\w,q}&=
 \frac{\pi ^6}{13824 h^6}-\frac{7 \pi ^4}{768 h^5}+\frac{21 \pi ^2}{128 h^4}+\frac{21}{64 h^3}+\mathrm{e.s.t.} \\
\langle\p_1^4\rangle_{\w,q}&=\frac{\pi ^8}{331776 h^8}-\frac{11 \pi ^6}{13824 h^7}+\frac{77 \pi ^4}{1536 h^6}-\frac{77 \pi ^2}{128 h^5}-\frac{231}{256 h^4}
+\mathrm{e.s.t.}\\
\langle\p_3\rangle_{\w,q}&=\frac{7 \pi ^4}{960 h^4}\\
\langle\p_1\p_3\rangle_{\w,q}&=\frac{7 \pi ^6}{23040 h^6}-\frac{7 \pi ^4}{256 h^5}+\mathrm{e.s.t.}\\
\langle\p_1^2\p_3\rangle_{\w,q}&=\frac{7 \pi ^8}{552960 h^8}-\frac{133 \pi ^6}{46080 h^7}+\frac{133 \pi ^4}{1024 h^6}+\mathrm{e.s.t.}\\
\langle\p_1^3\p_3\rangle_{\w,q}&=\frac{7 \pi ^{10}}{13271040 h^{10}}-\frac{161 \pi ^8}{737280 h^9}+\frac{3059 \pi ^6}{122880 h^8}-\frac{3059 \pi ^4}{4096 h^7}+\mathrm{e.s.t.}
\end{align*}
Here, ``e.s.t.'' stands for exponentially small terms.

Instead of \eqref{eq:needforp1}, we get
 \begin{multline*}
\left\langle \left(\p_1-\langle\p_1\rangle_{\w,q}\right)^4\right\rangle_{\w,q}=
 \frac{\pi ^4}{48 h^6}-\frac{7 \pi ^2}{8 h^5}-\frac{21}{16 h^4}+\cdots\\
 \neq
 \frac{\pi ^4}{48 h^6}+\frac{\pi ^2}{8 h^5}+\frac{3}{16 h^4}+\cdots
 =3\left\langle \left(\p_1-\langle\p_1\rangle_{\w,q}\right)^2\right\rangle_{\w,q}^2.
 \end{multline*}
Instead of \eqref{eq:needforp1p3}, we get
\begin{multline*}
 \left\langle \left(\p_1-\langle\p_1\rangle_{\w,q}\right)^3
\left(\p_3-\langle\p_3\rangle_{\w,q}\right)
\right\rangle_{\w,q}\\
=\frac{7 \pi ^6}{960 h^8}-\frac{273 \pi ^4}{320 h^7}+\cdots
\neq\frac{7 \pi ^6}{960 h^8}+\frac{7 \pi ^4}{320 h^7}+\cdots\\
=3\left\langle \left(\p_1-\langle\p_1\rangle_{\w,q}\right)^2\right\rangle_{\w,q}
\left\langle \left(\p_1-\langle\p_1\rangle_{\w,q}\right)\left(\p_3-\langle\p_3\rangle_{\w,q}\right)\right\rangle_{\w,q}.
\end{multline*}

From these sample computations, one would conclude that the convergence is not asymptotically gaussian and that no Wick-type theorem applies to this distribution.

\begin{remark}
In a previous version of this paper that was posted on the arXiv, and in the thesis \cite{mythesis}, there appeared different values for the expectations presented above. The reason for the change is that in the Spring of 2014 we discovered a defect in the computer program we were using. This is also the reason why we have added Appendix \ref{apx:program}, where we hope to make the whole procedure transparent.
\end{remark}

\section{Characters of near-involutions}
\label{sec:nearinvolutions}

Proposition \ref{prop:formulaforg} follows from Proposition \ref{prop:formula} below, and from the fact that
\[\g_\nu(\lambda)=\frac{2^{|\nu|/2}(|\lambda|/2)!}{\mathfrak z(\nu)\left(\frac{|\lambda|-|\nu|}2\right)!}
\frac{\chi^\lambda(\nu,2,2,\dots,2)}{\chi^\lambda(2,2,\dots,2)}.\]

\begin{definition}[Shifted Schur functions]\label{def:shiftedschur}
Let 
\[(x\downharpoonright k)=x(x-1)(x-2)\cdots(x-k+1).\]
The \emph{shifted Schur polynomials}
\[\s_\mu(x_1,\dots,x_n)=\frac{\det[(x_i+n-i\downharpoonright\mu_j+n-j)]}{\det[(x_i+n-i\downharpoonright n-j)]},\quad1\leq i,j\leq n\]
satisfy a stability condition
\[\s_\mu(x_1,\dots,x_n,0)=s_\mu(x_1,\dots,x_n),\]
which allows us to take inverse limits, just as in de definition of symmetric functions (see for example \cite{macdonald}). The resulting objects, $\s_\mu(x_1,x_2,\dots)$, are known as \emph{shifted Schur functions} \cite{okounkovolshanski} (or \emph{Frobenius-Schur functions} \cite{frobeniusschur}).
\end{definition}

\begin{proposition} \label{prop:formula}
For $\lambda$ and $\nu$ balanced partitions, with the Young diagram of $\nu$ entirely contained inside the Young diagram of $\lambda$,
\[\frac{\chi^{\lambda}(\nu,2,2,\dots,2)}{\chi^\lambda(2,2,\dots,2)}
=\frac{(|\lambda/\nu|/2)!}{(|\lambda|/2)!}
\sum_\mu(-1)^{o_\mu/2}\chi^\mu(\nu)\s_a(\alpha)\s_b(\beta),
\]
where the sum is taken over all balanced partitions $\mu$ with $|\mu|=|\nu|$ and whose Young diagrams are completely contained in the diagram of $\lambda$,
and the pairs of partitions $(\alpha,\beta)$ and $(a,b)$ are the 2-quotients of $\lambda$ and $\mu$, respectively.
\end{proposition}
\begin{remark}
This, together with the formula from the following lemma (in the case of $\mu$ empty) and together with any of the available expressions for the dimension $\dim \lambda$ of a partition (such as the Frame-Robinson-Thrall hook formula), gives an explicit
formula for $\chi^{\lambda}(\nu,2,2,\dots,2)$.
\end{remark}
\begin{lemma}\label{lem:dimension}
For $\lambda$ and $\mu$ balanced, and the Young diagram of $\mu$ completely contained inside the Young diagram of $\lambda$,
\[\chi^{\lambda/\mu}(2,2,\dots,2)=(-1)^{(o+o')/2}\binom{|\lambda/\mu|/2}{|\alpha/a|}
\dim(\alpha/a)\dim(\beta/b),\]
where $o$ and $o'$ are the numbers of odd parts of $\lambda$ and $\mu$, respectively, and the pairs of partitions $(\alpha,\beta)$ and $(a,b)$ are the 2-quotients of $\lambda$ and $\mu$, respectively.
\end{lemma}
\begin{proof}[Proof of Lemma \ref{lem:dimension}]
The is a straightforward consequence of the definition of the 2-quo\-tients and of the Murnaghan-Nakayama rule; see \cite{mythesis} for details.
\end{proof}

\begin{proof}[Proof of Proposition \ref{prop:formula}]
First, observe that the Murnaghan-Nakayama rule implies that
\[\chi^\lambda(\nu,2,2,\dots,2)=\sum_{|\mu|=|\nu|}\chi^{\mu}(\nu) \chi^{\lambda/\mu}(2,2,\dots,2).\]
The sum is over all partitions $\mu$ of size $|\nu|$ whose diagram is completely contained inside the diagram of $\lambda$. Clearly, $\chi^{\lambda/\mu}(2,2,\dots,2)$ vanishes unless $\mu$ is balanced, so all sums from this point on will be over balanced partitions $\mu$ of size $|\nu|$ and contained in $\lambda$. To this expression, we apply the formula from Lemma \ref{lem:dimension}, and we apply the same formula with empty $\mu$ to the denominator, to get
\[\frac{\chi^\lambda(\nu,2,2,\dots,2)}{\chi^\lambda(2,2,\dots,2)}=\sum_\mu (-1)^{o/2}B_\mu\chi^\mu(\nu)\frac{\dim(\alpha/a)}{\dim\alpha} \cdot\frac{\dim(\beta/b)}{\dim\beta}.
\]
where $o$ is the number of odd parts in $\mu$; $\alpha$, $\beta$, $a$, and $b$ are as in the statement of the proposition, and
\[B_\mu=\binom{|\lambda/\mu|/2}{|\alpha/a|}/\binom{|\lambda|/2}{|\alpha|}.\]
Now we apply, to each of the quotients of dimensions, the Okounkov-Olshanski formula \cite[Equation 0.14]{okounkovolshanski}
\[\frac{\dim(\lambda/\mu)}{\dim\lambda}
=\frac{\s_\mu(\lambda)}{n(n-1)\cdots(n-k+1)}\]
where $n=|\lambda|$ and $k=|\mu|$. We get exactly the formula in the proposition because
\[\frac{|\lambda|}2-|\alpha|=|\beta|\qquad\textrm{and}\qquad \frac{|\lambda/\mu|}2-|\alpha/a|=|\beta/b|.\qedhere\]
\end{proof}

\section{Moduli spaces}
\label{sec:application}
The function
\[\g_\nu(\lambda)=\frac{\f_{(\nu,2,2,\dots)}(\lambda)}{\f_{(2,2,\dots)}(\lambda)}\]
appears in the work of Eskin-Okounkov \cite{pillow} as an essential ingredient for the computation of the volumes of the different strata of the moduli space of quadratic differentials. The partition $\nu$ determines the corresponding stratum: the multiplicities of the zeros of the quadratic differentials are encoded as $\nu_i-2$. Simple poles are allowed and correspond to parts $\nu_i=1$.

It was shown in \cite{pillow} that these volumes are equal to the first term of the asymptotics of the expectations
\[\langle \g_\nu\rangle_{\w,q}\]
as $q\to1$. See also \cite{mythesis} for a detailed account.

\begin{remark}\label{rmk:quasimodularity}
In Eskin-Okounkov \cite{pillow}, it is proved that these expectations are quasimodular forms. Using Proposition \ref{prop:formulaforg}, we break up these expectations into linear combinations of expectations of the form
\begin{equation}\label{eq:expectationofschur}
\langle \s_a(\alpha)\s_b(\beta)\rangle_{\w,q}.
\end{equation}
 Using methods similar to those of \cite[Section 3]{pillow}, it can be proven \cite[Section 2.5]{mythesis} that these are also quasimodular when taken individually. This arguably shortens the proof of the quasimodularity of $Z_\nu(q)$ since most of Section 2 of \cite{pillow} becomes superfluous.
\end{remark}

\begin{remark}\label{rmk:matrix}
Recall \cite[Section 1.7]{macdonald} that the matrix $M$ relating the symmetric power functions $p_\nu$ to the Schur functions $s_\mu$ (the classical ones, not the shifted ones) with $|\nu|=|\mu|$ is the character table of the symmetric group $S(n)$, $M=\left(\chi^\mu(\nu)\right)_{\mu,\nu}$,
\[p_\nu=\sum_\mu\chi^\mu(\nu)s_\mu.\]
Its inverse is derived from
\[s_\mu=\sum_\nu \frac{\chi^\mu(\nu)}{\mathfrak z(\nu)}p_\nu,\]
to be the matrix $M^{-1}=\left(\chi^\mu(\nu)/\mathfrak z(\nu)\right)_{\mu,\nu}$.

If we define $V_\nu=2^{-|\nu|/2}\langle\g_\nu\rangle_{\w,q}$ and $H_\mu=\langle \s_a(\alpha)\s_b(\beta)\rangle_{\w,q}$, it follows from the formula in Proposition \ref{prop:formulaforg} that (asymptotically)
\[V_\nu=\sum_\mu \frac{\chi^\mu(\nu)}{\mathfrak z(\nu)} H_\mu. \]
In other words, the matrix relating $V_\nu$ to $H_\mu$ is the transpose of the inverse, $M^{-T}$. Whence we also know that (asymptotically)
\[H_\mu=\sum_\nu \chi^\mu(\nu)V_\nu.\]
In this way, Proposition \ref{prop:formulaforg} reveals some structure in the problem of determining $\langle\g_\nu\rangle_{\w,q}$ in general.
\end{remark}

We shall now derive some results about the asymptotic analysis of $\langle\g_\nu\rangle_{\w,q}$.

Let $C_1$ be the space of functions $r:\R\to\R$, $r(x)\geq|x|$, with Lipschitz constant $\leq 1$, that is, for all $x,y\in\R$,
\[|r(x)-r(y)|\leq |x-y|,\]
and such that
\[\int_\R[r(x)-|x|]\,dx=1.\]
Note that all contours $L_\lambda$ of partitions $\lambda$ belong to $C_1$, as do all possible limit shapes. We endow $C_1$ with the topology of the supremum norm.

For a function $f:C_1\to\R$ and a partition $\lambda$, we define $f(\lambda)$ to be the function evaluated on the contour of the partition, $f(L_\lambda)$.

Fixing $n$, let
\[\langle f\rangle_{\w}=\frac{1}Z\sum_{|\lambda|=n}\w(\lambda)f(\lambda),\]
where $Z=\sum_{|\lambda|=n}\w(\lambda)$.

\begin{corollary}\label{coro:expectations}
Let $f,g:C_1\to\R$ be two functions with finite $\w$-ex\-pec\-ta\-tions at each level $n$. Assume that $|f(\lambda)|$ and $|g(\lambda)|$ increase at most polynomially as $|\lambda|\to\infty$. In particular for $g$,
\[g(\lambda)=O(|\lambda|^b) \quad as \quad|\lambda|\to \infty\]
for some $b>0$. Assume additionally that $g(\lambda)|\lambda|^{-b}$ is continuous with the topology of $C_1$. Then
\[\lim_{n\to\infty}\frac1{n^b}\left(\langle f(\alpha)g(\beta)\rangle_\w-\langle f(\alpha)g(\alpha)\rangle_\w\right)=0,\]
where $(\alpha,\beta)$ is the 2-quotient of the partition $\lambda$ over which the sum of the expectation is taken.
\end{corollary}
\begin{proof}[Proof of Corollary \ref{coro:expectations}]
Let $\varepsilon>0$. Let $\delta>0$ be such that if $\|L_\alpha-L_\beta\|<\delta$ then $|g(L_\alpha)-g(L_\beta)|<\varepsilon n^b$. We want to show that the following tends to 0:
\[\frac{1}{Zn^b}\sum_{|\lambda|=n}\w(\lambda)f(\alpha)(g(\beta)-g(\alpha)).\]
Here, the summation can be split into two parts,
\[\sum_{\|L_\alpha-L_\beta\|<\delta}+\sum_{\|L_\alpha-L_\beta\|\geq\delta}.\]
The later can be bounded easily using Proposition \ref{prop:closetodiagonal}: we get
\[\frac{1}{Zn^b}\sum_{\|L_\alpha-L_\beta\|\geq\delta}O\left(e^{-K_\delta\sqrt n}\right) f(\alpha)(g(\beta)-g(\alpha)),\]
for some constant $K_\delta$ that depends on $\delta$.
Since $f$ and $g$ grow polynomially, this tends to 0 as $n\to\infty$. On the other hand, the first part in the summation above goes like
\[\frac{1}{Zn^b}\sum_{\|L_\alpha-L_\beta\|<\delta} O\left(e^{-K_\delta\sqrt n}\right) f(\alpha)\varepsilon n^b.\qedhere\]
\end{proof}

\begin{remark}\label{rmk:asymptoticsoftwoquots}
It follows from Corollary \ref{coro:expectations} and Proposition \ref{prop:limitshape} that the components $\alpha$ and $\beta$ of the 2-quotient of $\lambda$ also have a limit shape, which in fact coincides with the one corresponding to the uniform distribution.
\end{remark}

\begin{proposition}\label{prop:application}
With the same notations as above
\begin{align}
\label{eq:term1}\langle \g_\nu\rangle_{\w}&=\frac{K^{|\nu|}}{\mathfrak z(\nu)} \sum_{\mu} (-1)^{o_\mu/2} \chi^\mu(\nu)\frac{\dim a}{|a|!}\frac{\dim b}{|b|!}n^{n/2}\\
\label{eq:term2}&\quad+\frac{2K^{|\nu|-2}\langle\p_2\rangle_\w}{\mathfrak z(\nu)}
\sum_{\mu} (-1)^{o_\mu/2}\chi^\mu(\nu)v(a)v(b)\\
&\notag\quad+\textrm{terms of lower degree,}
\end{align}
where $K>0$, both sums are taken over balanced partitions $\mu$ of size $|\mu|=|\nu|$ and
\[v(\eta)=\frac{\dim\eta}{(|\eta|-2)!}\frac{\s_{(2)}(\eta)-\s_{(1,1)}
(\eta)}{|\eta|!}.\]
Moreover, the first term \eqref{eq:term1} of $\langle\g_\nu\rangle_\w$ always vanishes.
\end{proposition}
\begin{proof}[Proof of Proposition \ref{prop:application}]
First we will show that the term \eqref{eq:term1} vanishes. Since, for $\mu$ balanced \cite[Section 2.1]{mythesis},
\[\chi^\mu(2,2,\dots,2)=(-1)^{o_\mu}\binom{|\mu|}{|a|}\dim a\dim b\]
and $|a|+|b|=|\mu|/2$, the first term \eqref{eq:term1} becomes
\[\frac{K^{|\nu|}}{\mathfrak z(\nu)(|\nu|/2)!}\sum_\mu\chi^\mu(\nu)\chi^\mu(2,2,\dots,2),\]
which vanishes because of character orthogonality.

To prove that the highest degree term equals \eqref{eq:term1}, we start with the formula from Proposition \ref{prop:formulaforg}, we apply Corollaries \ref{coro:multiplicative} and \ref{coro:expectations},  and we recall from \cite{okounkovolshanski} that\footnote{See also \cite{frobeniusschur} for the full expansion of $\s_\mu$ in terms of the regular Schur functions $s_\lambda$.}
\[\s_\eta=\sum_\rho \frac{\chi^\eta(\rho)}{\mathfrak z(\rho)}\p_\rho+\textrm{lower order terms,}\]
where the sum is taken over all partitions $\rho$ of size $|\eta|$.
We have the asymptotics \cite[Section 3.3.5]{pillow}
\[\langle \p_\mu\rangle_\w=O\left(n^{\frac{|\mu|+\ell(\mu)}2}\right)\quad \textrm{as $n\to\infty$},\]
and a similar version for the expectation $\langle \p_\mu(\alpha)\rangle_\w$ can be obtained easily using a similar argument with the generating function in \cite[Section 2.5]{mythesis}.
This then implies that the terms that grow fastest are precisely those involving $\p_{(1,1,\dots,1)}$ only. This means that the highest degree term is
\[\frac{2^{|\nu|/2}}{\mathfrak z(\nu)}\sum_\mu(-1)^{o_\mu/2}\chi^\mu(\nu) \frac{\chi^a(1,\dots,1)}{\mathfrak z(1,\dots,1)}
\frac{\chi^b(1,\dots,1)}{\mathfrak z(1,\dots,1)}\langle \p_{(1,\dots,1)}(\alpha)\rangle_\w^2,
\]
which is equal to \eqref{eq:term1}.

The proof that the second term \eqref{eq:term2} is similar, taking now the terms involving only $\p_{(2,1,1,\dots\,1)}$ and using additionally the fact that
\begin{align*}
\chi^\eta(2,1,1,\dots,1)&=\chi^{(2)}(2)\dim(\eta/(2))+\chi^{(1,1)}(2)\dim(\eta/(1,1)) \\
&=\frac{\dim\eta}{|\eta|!/2}\left(\s_{(2)}(\eta)-\s_{(1,1)}(\eta)\right),
\end{align*}
which follows from the Murnaghan-Nakayama rule, and from the O\-koun\-kov-Ol\-shan\-ski \cite{okounkovolshanski} formula,
\[\frac{\dim(\lambda/\mu)}{\dim \lambda}=\frac{|\mu|!} {|\lambda|!}\s_\mu(\lambda).\qedhere\]
\end{proof}

\begin{remark}
It is possible to prove Proposition \ref{prop:application} without using Corollary \ref{coro:multiplicative}; see \cite[Sections 3.5 and 3.6]{mythesis}.
\end{remark}

\appendix
\section{Program listings}
\label{apx:program}
In this section we present a program that we have used to compute the asymptotic behavior of the expectations presented in Section \ref{sec:noclt}. The program is a script designed to run in Sage \cite{sage}.

What the program does is essentially a linear regression. We know from the work of Okounkov and Eskin \cite{pillow} that $\langle f\rangle_{\w,q}$
is a quasimodular form. Here, $f$ is the desired function and it must be in the algebra generated by the shifted power functions $\p_n$; $f$ appears in the program below with the name \vari{function\_of\_interest}. We compute the first few coefficients of this series, as well as the first few elements of a basis of the space of quasimodular forms, and then we invert the corresponding matrix to find out what the linear combination is. Our basis is given by products of the Eisenstein series $E_2(q^2),E_2(q^4),E_4(q^4)$, and it is hence very easy to compute the asymptotic behavior of each term.

The hardest part is the computation of the coefficients of the series $Z\langle f\rangle_{\w,q}$ because each of them is a combinatorial sum over all pairs of partitions whose sizes add up to the corresponding power of $q$. This is done by the slave program, which we present in Section \ref{sec:slave}. The computation of the basis and the linear regression are done by the main program, presented in Section \ref{sec:mainprogram}.

The output of the program is a polynomial in the variable \vari{h}. In the notations of Section \ref{sec:noclt}, \vari{h} is really a placeholder for $h^{-1}$. This polynomial thus respresents the non-exponentially small part the asymptotic behavior of $\langle f\rangle_{\w,q}$. If too few coefficients of $\langle f\rangle_{\w,q}$ have been computed with the slave program, the main program will effectively return a lot of trash.

\subsection{Slave program}
\label{sec:slave}
The following function computes very quickly the value of $\w(\lambda)$, and it takes three partitions: \vari{p} equal to $\lambda$, and also \vari{p1} and \vari{p2} equal to $\alpha$ and $\beta$, the components of the 2-quotient of $\lambda$.
\begin{lstlisting}
def pillowcaseweight(p,p1,p2): 
    a = p.hook_product(1)
    b = p1.hook_product(1)
    c = p2.hook_product(1)
    return a^2/(b*c)^4/2^(2*sum(p))
\end{lstlisting}
The following function computes $\p_k(\lambda)$. It takes, as single argument, the partitions \vari{l} equal to $\lambda$. The variable \vari{z} is a dictionary with the values of $\zeta(-k)$.
\begin{lstlisting}
def shiftedpower(k,l): 
    z={1:-1/12,2:0,3:1/120,4:0,5:-1/252,6:0,7:1/240}
    box = QQ((1-2**(-k))*z[k])
    for i in range(len(l)):
        box += (l[i]-i-1/2)^k-(-i-1/2)^k
    return box
\end{lstlisting}
The following function computes the value of 
\[\sum_{|\alpha|=n_1,|\beta|=n_2} \w(\lambda)f(\lambda)\]
for fixed integers $n_1$ and $n_2$. Here, $\lambda$ is the partition with 2-quotient $(\alpha,\beta)$, and $f$ is a function (\vari{function\_of\_interest}, below) that must be defined separately. Then, it saves it to a file. The function takes two positive integers \vari{n1} and \vari{n2}, and a string, \vari{filename}, that determines the name of the file where the result is saved. The function does nothing if the file already exists.
\begin{lstlisting}
def work(n1,n2,filename):
    if not os.path.exists(filename):
        print('  working...')
        box=0
        for p1 in Partitions(n1):
            for p2 in Partitions(n2):
                p=Partition( \
                    core=[],quotient=[p1,p2])
                pval=function_of_interest(p1,p2)
                theweight=pillowcaseweight(p,p1,p2)
                box+=theweight*pval
        with open(filename,'w') as f:
            f.write(str(box)+"\n")
            f.close()
        print('     done!')
    else:
        print('  was already done')
\end{lstlisting}
The following function looks for files named ``$n_1$-$n_2$.dat'' with $n_1+n_2$ less than or equal to the parameter \vari{n}, and it adds up their contents to obtain the value of
\[\sum_{|\lambda|=n_1+n_2}\w(\lambda)f(\lambda).\]
It stores a list of the results for each different value of $n_1+n_2$ in a file specified by the argument \vari{filename}.
\begin{lstlisting}
def make_cache(n,filename):
    rows=[]
    for i in range(0,n+1):
        row=0
        for k in range(i+1):
            print(k,i-k)
            with open(str(k)+'-'+ \ 
                    str(i-k)+'.dat','r') as dat:
                datum=dat.read()
                dat.close()
                print(datum)
                val=sage_eval(datum)
                row+=val
        rows.append(row)
    with open(filename,'w') as f:
        for r in rows:
            f.write(str(r)+'\n')
        f.close()
\end{lstlisting}
The following is the function one would call to get the whole job done, as it automatizes the whole proceeding of the slave program. Given the integer \vari{n}, it produces a file with the list of $\w$-expectations of $f$ over the partitions of size $|\lambda|=m$ for $m\leq n$. This file is what the main program needs in order to finish the computation. 
\begin{lstlisting}
def bulkwork(n,filename):
    for m in range(n+1):
        for i in range(m+1):
            print(str(i)+'-'+str(m-i))
            work(i,m-i,str(i)+'-'+str(m-i)+'.dat')
        print('Done '+str(m))
    make_cache(n,filename)
\end{lstlisting}

\subsection{Main program}
\label{sec:mainprogram}
\begin{lstlisting}

# It is necessary to specify the name of the file
# generated by the slave.
filename="filename" 
 
T=PolynomialRing(QQ,'q')
q=T.gen()
S=PolynomialRing(QQ,'h')
h=S.gen()

def Eisenstein(k,qq):  
    # This is E_2k(q); good for k=1,2 only!
    constantterm={1:-1/24,2:1/240}
    return constantterm[k] + \
        sum([sigma(i+1,2*k-1)*qq^(i+1) \
        for i in range(maxdegree)])

# Load the first few coefficients of the expectation
# of interest from the file, and turn them into a 
# polynomial.
poly=0 
with open(filename,'r+') as f:
    readin=1
    counter=0
    readinstr=" "
    while len(readinstr) > 0:
        readinstr=f.readline()
        if len(readinstr)>0:
             readin = QQ(readinstr)
             poly+=q^(2*counter)*readin
             counter+=1
maxdegree=2*counter-2
print("Maximum degree found: "+str(maxdegree))
    
# Form the inverse of the partition function Z.
sqrrt=sum([(-q)^n*binomial(1/2,n) \
    for n in range(maxdegree+1)])
inverseZ= \
    prod([sqrrt(q^(2*i)).truncate(maxdegree+1) \
    for i in range(1,maxdegree+1)])

# List of all the different products of 
# Eisenstein series.
Elist=[T(1)] # Accumulates products; 
             # starts with the empty product.
products=[[]] # Keeps track of what each product
              # is made of.
factor={1:Eisenstein(1,q^2),2:Eisenstein(1,q^4), \
    3:Eisenstein(2,q^4)}           
for i in range(1,maxdegree/2+1):
    for p in Partitions(i,max_part=3):
        products.append(p)
        Elist.append(prod([factor[part] for part \
            in p]).truncate(maxdegree+1))
# Extract their coefficients and save them 
# in a matrix
coeffmatrix=[]
for pol in Elist:
    pollist=pol.list()
    pollist.extend([0 for i in \
        range(maxdegree+1-len(pollist))])
    coeffmatrix.append(pollist[2:len(pollist)+1:2]) 
cfmatrix=matrix(QQ,coeffmatrix)
# Extract a linearly-independent subset.
independentEs=cfmatrix.pivot_rows()

# Do linear regression.
thematrix= \ 
    (cfmatrix.matrix_from_rows( \
        independentEs)).transpose()
prevector=(poly*inverseZ).list()
x=thematrix.inverse()*vector(QQ, \ 
    prevector[2:maxdegree+1:2])
# Print out the linear decomposition of the
# generating function into the products of 
# Eisenstein series.
for i in range(len(x)):
    print(str(x[i])+" "+ \
        str(products[independentEs[i]]))

# Find the asymptotic value.
# NB: h here is h^(-1) in the text, 
# i.e. q=exp(h) here.
asymptoticvalue = {1: pi^2 / 24 * h + h / 4, \
    2: pi^2 / 96 * h^2 + h / 8, \
    3: pi^4 / 3840 * h^4}
    # These are the asymptotic values of E_2(q^2), 
    # E_2(q^4), E_4(q^4), respectively.
for i in range(len(x)):
    thissummand=x[i]
    for p in products[independentEs[i]]:
        thissummand*=asymptoticvalue[p]
    box+=thissummand
print "The asymptotic value is " + str(expand(box))
\end{lstlisting}



\def\cprime{$'$}
\providecommand{\bysame}{\leavevmode\hbox to3em{\hrulefill}\thinspace}
\providecommand{\MR}{\relax\ifhmode\unskip\space\fi MR }
\providecommand{\MRhref}[2]{%
  \href{http://www.ams.org/mathscinet-getitem?mr=#1}{#2}
}
\providecommand{\href}[2]{#2}
\begin{thebibliography}{10}

\bibitem{andrews}
George~E. Andrews, \emph{The theory of partitions}, Addison-Wesley Publishing
  Co., Reading, Mass.-London-Amsterdam, 1976, Encyclopedia of Mathematics and
  its Applications, Vol. 2. \MR{0557013 (58 \#27738)}

\bibitem{rightangledbilliards}
Jayadev~S. Athreya, Alex Eskin, and Anton Zorich, \emph{Right-angled billiards
  and volumes of moduli spaces of quadratic differentials on
  {${\mathbb{C}P^1}$}}, Appendix by Jon Chaika. arXiv:1212.1660.

\bibitem{jenkinsstrebeldifferentials}
\bysame, \emph{Counting generalized {J}enkins-{S}trebel differentials}, Geom.
  Dedicata \textbf{170} (2014), 195--217. \MR{3199485}

\bibitem{billingsley}
Patrick Billingsley, \emph{Probability and measure}, anniversary ed., Wiley
  Series in Probability and Statistics, John Wiley \& Sons Inc., Hoboken, NJ,
  2012, With a foreword by Steve Lalley. \MR{2893652}

\bibitem{blochokounkov}
Spencer Bloch and Andrei Okounkov, \emph{The character of the infinite wedge
  representation}, Adv. Math. \textbf{149} (2000), no.~1, 1--60. \MR{1742353
  (2001g:11059)}

\bibitem{pillow}
A.~Eskin and A.~Okounkov, \emph{Pillowcases and quasimodular forms}, Algebraic
  geometry and number theory, Progr. Math., vol. 253, Birkh\"auser Boston,
  Boston, MA, 2006, pp.~1--25. \MR{2263191 (2007g:14065)}

\bibitem{eskinmasurzorich}
Alex Eskin, Howard Masur, and Anton Zorich, \emph{Moduli spaces of abelian
  differentials: the principal boundary, counting problems, and the
  {S}iegel-{V}eech constants}, Publ. Math. Inst. Hautes \'Etudes Sci. (2003),
  no.~97, 61--179. \MR{2010740 (2005b:32029)}

\bibitem{branchedcoverings}
Alex Eskin and Andrei Okounkov, \emph{Asymptotics of numbers of branched
  coverings of a torus and volumes of moduli spaces of holomorphic
  differentials}, Invent. Math. \textbf{145} (2001), no.~1, 59--103.
  \MR{1839286 (2002g:32018)}

\bibitem{eskinokounkovpandharipande}
Alex Eskin, Andrei Okounkov, and Rahul Pandharipande, \emph{The theta
  characteristic of a branched covering}, Adv. Math. \textbf{217} (2008),
  no.~3, 873--888. \MR{2383889 (2008k:14065)}

\bibitem{ivanovolshanski}
Vladimir Ivanov and Grigori Olshanski, \emph{Kerov's central limit theorem for
  the {P}lancherel measure on {Y}oung diagrams}, Symmetric functions 2001:
  surveys of developments and perspectives, NATO Sci. Ser. II Math. Phys.
  Chem., vol.~74, Kluwer Acad. Publ., Dordrecht, 2002, pp.~93--151. \MR{2059361
  (2005d:05148)}

\bibitem{kontsevichlyapunovexponents}
M.~Kontsevich, \emph{Lyapunov exponents and {H}odge theory}, The mathematical
  beauty of physics ({S}aclay, 1996), Adv. Ser. Math. Phys., vol.~24, World
  Sci. Publ., River Edge, NJ, 1997, pp.~318--332. \MR{1490861 (99b:58147)}

\bibitem{kontsevich}
Maxim Kontsevich, \emph{Intersection theory on the moduli space of curves and
  the matrix {A}iry function}, Comm. Math. Phys. \textbf{147} (1992), no.~1,
  1--23. \MR{1171758 (93e:32027)}

\bibitem{kontsevichzorich}
Maxim Kontsevich and Anton Zorich, \emph{Connected components of the moduli
  spaces of {A}belian differentials with prescribed singularities}, Invent.
  Math. \textbf{153} (2003), no.~3, 631--678. \MR{2000471 (2005b:32030)}

\bibitem{macdonald}
I.~G. Macdonald, \emph{Symmetric functions and {H}all polynomials}, second ed.,
  Oxford Mathematical Monographs, The Clarendon Press Oxford University Press,
  New York, 1995, With contributions by A. Zelevinsky, Oxford Science
  Publications. \MR{1354144 (96h:05207)}

\bibitem{masur}
Howard Masur, \emph{Interval exchange transformations and measured foliations},
  Ann. of Math. (2) \textbf{115} (1982), no.~1, 169--200. \MR{644018
  (83e:28012)}

\bibitem{mumford}
D.~Mumford, \emph{Tata lectures on theta. {I}}, Modern Birkh\"auser Classics,
  Birkh\"auser Boston Inc., Boston, MA, 2007, With the collaboration of C.
  Musili, M. Nori, E. Previato and M. Stillman, Reprint of the 1983 edition.
  \MR{2352717 (2008h:14042)}

\bibitem{gromovwittenhurwitz}
A.~Okounkov and R.~Pandharipande, \emph{Gromov-{W}itten theory, {H}urwitz
  theory, and completed cycles}, Ann. of Math. (2) \textbf{163} (2006), no.~2,
  517--560. \MR{2199225 (2007b:14123)}

\bibitem{okounkovolshanski}
A.~Okun{\cprime}kov and G.~Ol{\cprime}shanski{\u\i}, \emph{Shifted {S}chur
  functions}, Algebra i Analiz \textbf{9} (1997), no.~2, 73--146. \MR{1468548
  (99f:05118)}

\bibitem{frobeniusschur}
Grigori Olshanski, Amitai Regev, and Anatoly Vershik, \emph{Frobenius-{S}chur
  functions}, Studies in memory of {I}ssai {S}chur ({C}hevaleret/{R}ehovot,
  2000), Progr. Math., vol. 210, Birkh\"auser Boston, Boston, MA, 2003, With an
  appendix by Vladimir Ivanov, pp.~251--299. \MR{1985729 (2004f:05188)}

\bibitem{mythesis}
Rodolfo R\'ios-Zertuche, \emph{Near-involutions, the pillowcase distribution
  and quadratic differentials}, Ph.D. thesis, Princeton University, 2012.

\bibitem{stanley}
Richard~P. Stanley, \emph{Enumerative combinatorics. {V}ol. 2}, Cambridge
  Studies in Advanced Mathematics, vol.~62, Cambridge University Press,
  Cambridge, 1999, With a foreword by Gian-Carlo Rota and appendix 1 by Sergey
  Fomin. \MR{1676282 (2000k:05026)}

\bibitem{sage}
W.\thinspace{}A. Stein et~al., \emph{{S}age {M}athematics {S}oftware ({V}ersion
  6.1.1)}, The Sage Development Team, 2014, {\tt http://www.sagemath.org}.

\bibitem{veech}
William~A. Veech, \emph{Moduli spaces of quadratic differentials}, J. Analyse
  Math. \textbf{55} (1990), 117--171. \MR{1094714 (92e:32014)}

\bibitem{vershikkerov1985}
A.~M. Vershik and S.~V. Kerov, \emph{Asymptotic behavior of the maximum and
  generic dimensions of irreducible representations of the symmetric group},
  Funktsional. Anal. i Prilozhen. \textbf{19} (1985), no.~1, 25--36, 96.
  \MR{783703 (86k:11051)}

\bibitem{wilf}
Herbert~S. Wilf, \emph{Generatingfunctionology}, third ed., A K Peters Ltd.,
  Wellesley, MA, 2006. \MR{2172781 (2006i:05014)}

\bibitem{zorich}
Anton Zorich, \emph{Flat surfaces}, Frontiers in number theory, physics, and
  geometry. {I}, Springer, Berlin, 2006, pp.~437--583. \MR{2261104
  (2007i:37070)}

\end{thebibliography}


\def\cprime{$'$}
\providecommand{\bysame}{\leavevmode\hbox to3em{\hrulefill}\thinspace}
\providecommand{\MR}{\relax\ifhmode\unskip\space\fi MR }
\providecommand{\MRhref}[2]{%
  \href{http://www.ams.org/mathscinet-getitem?mr=#1}{#2}
}
\providecommand{\href}[2]{#2}


\ACKNO{The author is deeply indebted to his doctoral adviser, Professor Andrei Okounkov, who suggested the problem and provided much valuable insight in the course of numerous conversations, and to Professor Grigori Olshanski for several useful discussions and for doing a full reading of the thesis this paper is based on.}


\end{document}